\documentclass{amsart}
\title[On supercongruences for truncated sums of squares...]
{On supercongruences for truncated sums of squares of basic hypergeometric series}
\usepackage{amssymb,amsmath,amsthm,epsfig,graphics,latexsym}
\usepackage{enumerate}
 \theoremstyle{definition}

  \theoremstyle{plain}
  \newtheorem{lemma}      {Lemma}
  
  \newtheorem{theorem}    {Theorem}
  
  \newtheorem{corollary}  {Corollary}
  
  \newtheorem{conjecture} {Conjecture}

  \newcommand{\fr}{\frac}

\begin{document}
  \author{Mohamed El Bachraoui}
  %\thanks{This research was supported by UAEU/COS grant IRG-17/15}
  %\address{Dept. Math. Sci,
 %United Arab Emirates University, PO Box 17551, Al-Ain, UAE}
\email{melbachraoui@uaeu.ac.ae}
\begin{abstract}
Congruences of truncated sums of infinite series do not directly extend to
congruences of the truncated sums of higher powers of these infinite series.
Guo and Zudilin recently established a variety of supercongruences for  truncated sums of certain basic
 hypergeometric series.
 In this note we extend some of these supercongruences to the truncated sums of the squares of the
 corresponding series.
\end{abstract}
\date{\textit{\today}}
 \maketitle
\section{Introduction}
\noindent
Recall that for a complex number $q$ and a complex variable $a$, the $q$-shifted factorials are given by
\[
(a;q)_0= 1,\quad (a;q)_n = \prod_{i=0}^{n-1}(1-a q^i),\quad
(a;q)_{\infty} = \lim_{n\to\infty}(a;q)_n =\prod_{i=0}^{\infty}(1-a q^i) \quad (|q|<1).
\]
The $q$-integer is given for any nonnegative integer $n$ by
\[
[n]_q=[n] = \fr{1-q^n}{1-q} = 1+q+\cdots+q^{n-1}.
%{n\brack m} = {n \brack m}_q =
%\frac{(q;q)_n}{(q;q)_m (q;q)_{n-m}} .
\]
Let
$\mathbb{Z}[q]$ denote the set of polynomials in $q$ with integer coefficients and let
$\mathbb{Z}(q)$ denote the set of rational functions in $q$ with integer coefficients.
The $n$-th cyclotomic polynomial is the polynomial in $\mathbb{Z}[q]$ given by
\[
\Phi_n(q) = \prod_{\substack{j=1 \\ \gcd(j,n)=1}}^n(q- \zeta^j) ,
\]
where $\zeta = e^{2\pi i/n}$ is the $n$-th root of unity.
Given polynomials $A_1(q), A_2(q), P(q) \in \mathbb{Z}[q]$, the congruence
$\fr{A_1(q)}{A_2(q)} \equiv 0 \pmod{P(q)}$ means that $P(q)$ divides $A_1(q)$ and that
$\gcd\big(P(q),A_2(q) \big) =1$; and in general for rational functions $A(q), B(q) \in\mathbb{Z}(q)$,
the congruence $A(q)\equiv B(q) \pmod{P(q)}$ means that $A(q)-B(q)\equiv 0\pmod{P(q)}$.
Many of the recent publications dealing with supercongruence relations were motivated by Ramanujan's formulas for $\fr{1}{\pi}$ as
in Ramanujan~\cite{Ramanujan} and Berndt~\cite{Berndt} and the later work by
Van~Hamme~\cite{VanHamme} where he
stated a variety of conjectures on the $p$-adic analogues of Ramanujan-like formulas for $1/\pi$. See for instance~\cite{Guillera-Zudilin, Guo-1, Guo-2, Guo-3, Guo-4, Osburn-Zudilin}.
For other related congruences, we refer to \cite{Guo-Wang, Shi-Pan, Straub, Zudilin}.
Guo~and~Zudilin~\cite[Theorems 1.1-1.4]{Guo-Zudilin-2} recently used the so-called \emph{$q$-microscope} approach to establish the following supercongruences.

\noindent
If $n$ is a positive odd integer, then
\begin{equation}\label{Eq:GZ-1}
\sum_{k=0}^{n-1}(-1)^k \fr{(q;q^2)_k (-q;q^2)_k^2}{(q^4;q^4)_k (-q^4;q^4)_k^2} [6k+1] q^{3k^2}
\equiv 0 \pmod{[n]}
\end{equation}
\noindent
and
\begin{equation}\label{Eq:GZ-2}
\sum_{k=0}^{n-1} \fr{(q^2;q^4)_k (-q;q^2)_k^2}{(q^4;q^4)_k (-q^4;q^4)_k^2} [6k+1] q^{k^2}
\equiv 0 \pmod{[n]},
\end{equation}
\noindent
and if $n$ is a positive integer which is coprime to $6$, and $a$ and $q$ are complex numbers, then
we have modulo $[n](1-aq^n)(a-q^n)$,
\begin{equation}\label{Eq:GZ-3}
\sum_{k=0}^{n-1}\fr{(aq;q^2)_k (q/a;q^2)_k (q;q^2)_{2k}}{(aq^6;q^6)_k (q^6/a;q^6)_k (q^2;q^2)_{2k}}
[8k+1] q^{2k^2}
\equiv q^{-(n-1)/2}[n]\Big(\fr{-3}{n}\Big),
\end{equation}
\noindent
where $\Big(\fr{-3}{n}\Big)$ is the Jacobi symbol, and as a consequence of (\ref{Eq:GZ-3}) by letting $a\to 1$
\begin{equation*}%\label{Eq:GZ-4}
\sum_{k=0}^{n-1} \fr{(q;q^2)_{k}^2 (q;q^2)_{2k}}{(q^6;q^6)_{k}^2 (q^2;q^2)_{2k}} [8k+1] q^{2k^2}
\equiv q^{-(n-1)/2}[n]\Big(\fr{-3}{n}\Big) \pmod{[n]\Phi_n(q)^2} .
\end{equation*}
An important step towards establishing their supercongruences,
Guo and Zudilin~\cite{Guo-Zudilin-1, Guo-Zudilin-2} evaluated the basic hypergeometric series corresponding
to the three congruences~(\ref{Eq:GZ-1})--(\ref{Eq:GZ-3}) respectively as follows:
\begin{equation}\label{Eq:GZ-1-1}
\sum_{k=0}^\infty(-1)^k \fr{(q;q^2)_k (-q;q^2)_k^2}{(q^4;q^4)_k (-q^4;q^4)_k^2}[6k+1] q^{3k^2}
=
\fr{(q^3;q^2)_\infty}{(-q^4;q^4)_\infty},
\end{equation}
\begin{equation}\label{Eq:GZ-2-1} %\cite{Guo-Zudilin-2}
\sum_{k=0}^\infty \fr{(q^2;q^4)_k (-q;q^2)_k^2}{(q^4;q^4)_k (-q^4;q^4)_k^2}[6k+1] q^{k^2}
=
\fr{(-q^2;q^4)_\infty^2}{(1-q)(-q^4;q^4)_\infty^2},
\end{equation}
and
\begin{equation}\label{Eq:GZ-3-1}
\sum_{k=0}^{\infty}\fr{(aq;q^2)_k (q/a;q^2)_k (q;q^2)_{2k}}{(aq^6;q^6)_k (q^6/a;q^6)_k (q^2;q^2)_{2k}}
[8k+1] q^{2k^2}
\end{equation}
\[
=
\fr{(q;q^2)_{\infty}(q^6;q^6)_{\infty} (aq^3;q^6)_{\infty} (q^3/a;q^6)_{\infty}}
{(1-q) (q^2;q^2)_{\infty}(q^3;q^6)_{\infty} (aq^6;q^6)_{\infty} (q^6/a;q^6)_{\infty}}.
\]
Roughly speaking, the $q$-microscope approach of Guo~and~Zudilin~\cite{Guo-Zudilin-2} goes as follows. To derive a supercongruence for a sum ranging from
$0$ to $n-1$, they first express the infinite version of their truncated sum
as an infinite product, say
\[ \sum_{k=0}^\infty c_q(k) = \prod_{k=0}^\infty C_k(q),
\]
and secondly they investigated the radial limit on both sides as $q$ approaches a primitive root
of unity $\zeta$ of degree a divisor $d$ of $n$.
Then after proving that
$\lim_{q\to \zeta}\prod_{k=0}^\infty C_k(q)$ is either zero or bounded above, they proceed with
other steps to obtain the desired information on $\sum_{k=0}^{n-1} c_{\zeta}(k)$.
However, this approach seems
not to extend directly to the truncated sums of higher powers of such series.
To clarify this difficulty, observe that while higher powers of the right-hand-sides of~(\ref{Eq:GZ-1-1})--(\ref{Eq:GZ-3-1}) behave nicely
as far as convergence when $q\to \zeta$ is concerned, the situation is quite unclear
when powers of the corresponding left-hand-sides are considered.
Our main goal in this note is to provide supercongruence relations for the truncated sums
ranging from $0$ to $n-1$ of the squares
of the series on the left-hand-sides of~(\ref{Eq:GZ-1-1})--(\ref{Eq:GZ-3-1}).
However, we were not able to obtain congruences for the sums ranging from $0$ to $\fr{n-1}{2}$. Among key ideas in the work of Guo~and~Zudilin~\cite{Guo-Zudilin-2} we find the following properties
of the $k$-th term $c_q(k)$ of the basic hypergeometric series which they studied:
\begin{equation}\label{Eq:key}
c_{\zeta}(k)=0 \ \text{for\ } \fr{d-1}{2} < k \leq d-1
\ \text{and \ }
\lim_{q\to\zeta} \fr{c_q(ld+k)}{c_q(ld)} = c_{\zeta}(k),
\end{equation}
where $d$ is a positive odd integer and $\zeta$ is a primitive $d$-th root of unity.
Our proofs rely on the $q$-microscope method~\cite{Guo-Zudilin-2} combined the two conditions in~(\ref{Eq:key})
and their consequences which we record in the following lemma.
\begin{lemma}\label{lem:elementary}
Let $d$ be a positive integer and let $\{c(k)\}_{k=0}^\infty$ be a sequence of complex numbers.

\noindent
(a)\ If  $c(k)=0$ for $\fr{d-1}{2} < k \leq d-1$, then
\begin{equation*}\label{Eq:lem-1}
\Big(\sum_{j=0}^{d-1} c(j) \Big)^2 = \sum_{k=0}^{d-1} \sum_{j=0}^k c(j) c(k-j).
\end{equation*}
\noindent
(b)\ If in addition, $\fr{c(ld+k)}{c(ld)} = c(k)$ for all nonnegetive integers $k$ and $l$ such that
$0\leq k < d-1$, then
\[
\sum_{j=0}^{ld+k} c(j) c(ld+k-j)
=
\sum_{i=0}^l c(i d)c\big( (l-i)d \big) \sum_{j=0}^k c(j)c(k-j).
\]
\end{lemma}
\noindent
Note that Lemma~\ref{lem:elementary} can be extended to the convolution of any two sequences which satisfy the conditions~(a) and~(b) in the lamma.
The rest of this paper is organized as follows. In Section~\ref{sec:results} we state our main theorem
and their $p$-adic consequences. In Sections~\ref{sec:thm1-proof}-\ref{sec:thm3-proof} we give
the proofs of Theorems~\ref{thm:main-1}-\ref{thm:main-3} respectively. Section~\ref{sec:lem-proof} is devoted
to the proof of Lemma~\ref{lem:elementary}. In Section~\ref{sec:conjecture} we give two conjectures on supercongruences which are related to our results.
\section{Statement of results}\label{sec:results}
\begin{theorem}\label{thm:main-1}
Let $c_q(k)$  be the $k$-th term of
\[
\sum_{k=0}^{\infty}(-1)^k \fr{(q;q^2)_k (-q;q^2)_k^2}{(q^4;q^4)_k (-q^4;q^4)_k^2} [6k+1] q^{3k^2}
\]
and  let 
\[
a_q(k) = \sum_{j=0}^k c_q(j) c_q(k-j).
\]
Then for any positive odd integer $n$ we have
\[
\sum_{k=0}^{n-1} a_q(k) \equiv 0 \pmod{[n]}.
\]
\end{theorem}
\noindent
By letting $q\to 1$ one obtains the following $p$-adic identity.
\begin{corollary}\label{cor:main-1}
For any odd prime number $p$, we have
\[
\sum_{k=0}^{p-1}\fr{(-1)^k}{8^{k}} \sum_{j=0}^k
{2j\choose j}{2k-2j\choose k-j}(6j+1)(6k-6j+1)
\equiv 0\pmod{p}.
\]
\end{corollary}
\begin{theorem}\label{thm:main-2}
Let $c_q(k)$  be the $k$-th term of
\[
\sum_{k=0}^{\infty}\fr{(q^2;q^4)_k (-q;q^2)_k^2}{(q^4;q^4)_k (-q^4;q^4)_k^2} [6k+1] q^{k^2}
\]
and  let 
\[
a_q(k) = \sum_{j=0}^k c_q(j) c_q(k-j).
\]
Then for any positive odd integer $n$ we have
\[
\sum_{k=0}^{n-1} a_q(k) \equiv 0 \pmod{[n]}.
\]
\end{theorem}
\noindent
By letting $q\to 1$ one obtains the following $p$-adic identity.
\begin{corollary}\label{cor:main-2}
For any odd prime number $p$, we have
\[
\sum_{k=0}^{p-1}\fr{1}{4^{k}} \sum_{j=0}^k
{2j\choose j}{2k-2j\choose k-j}(6j+1)(6k-6j+1)
\equiv 0\pmod{p}.
\]
\end{corollary}
\begin{theorem}\label{thm:main-3}
Let $c_q(k)$ be the $k$-th term
\[
 \sum_{k=0}^{\infty}\fr{(aq;q^2)_k (q/a;q^2)_k (q;q^2)_{2k}}{(aq^6;q^6)_k (q^6/a;q^6)_k (q^2;q^2)_{2k}}
[8k+1] q^{2k^2}
\]
and let \[
a_q(k) = \sum_{j=0}^k c_q(j) c_q(k-j).
\]
Then for any positive odd integer which is coprime to $6$ we have
\[
\sum_{k=0}^{n-1} a_q(k) \equiv q^{-(n-1)}[n]^2 \pmod{[n]\Phi_n(q)^2}.
\]
\end{theorem}
\noindent
Upon letting $a\to 1$ in Theorem~\ref{thm:main-3} we obtain the following special case.
\begin{theorem}\label{thm:main-4}
Let $c_q(k)$ be the $k$-th term of
\[
\sum_{k=0}^{\infty} \fr{(q;q^2)_{k}^2 (q;q^2)_{2k}}{(q^6;q^6)_{k}^2 (q^2;q^2)_{2k}} [8k+1] q^{2k^2}
\]
and let \[
a_q(k) = \sum_{j=0}^k c_q(j) c_q(k-j).
\]
Then for any positive odd integer which is coprime to $6$ we have
\[
\sum_{k=0}^{n-1} a_q(k) \equiv q^{-(n-1)}[n]^2 \pmod{[n]\Phi_n(q)^2}.
\]
\end{theorem}
\begin{corollary}\label{cor:main-4}
Let $p>3$ be a prime number. Then we have
\[
\sum_{k=0}^{p-1}\fr{1}{2^{8k} 3^{2k}} \sum_{j=0}^k
{2j\choose j}^2{4j\choose 2j} {2k-2j\choose k-j}^2{4k-4j\choose 2k-2j}(8j+1)(8k-8j+1)
\]
\[
\equiv p^2  \pmod{p^3}.
\]
\end{corollary}
\section{Proof of Theorem~\ref{thm:main-1}}\label{sec:thm1-proof}
\noindent
The result is evident for $n=1$. Now let $n>1$ be odd, let $d$ be an odd divisor of $n$, and let
$\zeta$ be a primitive root of unity of degree $d$. With the help of the following
basic properties
\begin{equation}\label{Eq:zeta}
(-\zeta;\zeta)_d = (-\zeta^2;\zeta^2)_d = (-\zeta;\zeta^2)_d = (-\zeta^4;\zeta^4)_d = 2,
\end{equation}
Guo~and~Zudilin~\cite{Guo-Zudilin-2} showed that for any nonnegative integer $l$
\begin{equation}\label{Eq:c-ld-1}
\lim_{q\to\zeta}c_q(ld) = c_{\zeta}(ld) = \fr{(-1)^l}{8^l}{2l\choose l} \quad\text{and so\quad}
\sum_{l=0}^\infty c_{\zeta}(ld) = \fr{\sqrt{6}}{3}.
\end{equation}
By the relation~(\ref{Eq:GZ-1-1}) we have
\[
\sum_{k=0}^\infty a_q(k) =\fr{(q^3;q^2)_\infty^2}{(-q^4;q^4)_\infty^2}.
\]
Note that for the right-hand-side of the previous identity, it is to see that
\begin{equation}\label{left-limit}
\lim_{q\to \zeta} \fr{(q^3;q^2)_\infty^2}{(-q^4;q^4)_{\infty}^2}
= \lim_{q\to \zeta} \Big(\fr{(q^3;q^2)_\infty}{(-q^4;q^4)_{\infty}} \Big)^2 = 0.
\end{equation}
As to the left-hand-side, we have
\[
\sum_{k=0}^{\infty} a_q(k) = \sum_{l=0}^{\infty} \sum_{k=0}^{d-1} a_q(ld+k)
= \sum_{l=0}^{\infty} \sum_{k=0}^{d-1} \sum_{j=0}^{ld+k} c_q(j) c_q(ld+k-j).
\]
\noindent
Now by (\ref{Eq:key}) the sequence $\{c_{\zeta}(k))_{k=0}^{\infty}$ satisfies the conditions
of Lemma~\ref{lem:elementary}(a, b) and so,
the limit as $q\to \zeta$ of the rightmost sum of the foregoing formula becomes
\[
\sum_{l=0}^{\infty} \sum_{k=0}^{d-1} \sum_{i=0}^l c_{\zeta}(id) c_{\zeta} \big( (l-i)d \big)
\sum_{j=0}^k c_{\zeta}(j)c_{\zeta}(k-j)
\]
\[
=
\sum_{l=0}^{\infty} \Big(\sum_{i=0}^l c_{\zeta}(id) c_{\zeta} \big( (l-i)d \big) \Big)
\sum_{k=0}^{d-1} a_{\zeta}(k)
=
\Big( \sum_{l=0}^{\infty} \fr{(-1)^l}{8^l} {2l\choose l} \Big)^2 \sum_{k=0}^{d-1} a_{\zeta}(k)
\]
\[
= \fr{6}{9} \sum_{k=0}^{d-1} a_{\zeta}(k)
= 0,
\]
where the second and third identities follow by~(\ref{Eq:c-ld-1}) and the last identity
follows from~(\ref{left-limit}). Thus we have
\[
\sum_{k=0}^{d-1} a_{\zeta}(k) = 0.
\]
Furthermore, we similarly get
\[
\begin{split}
\sum_{k=0}^{n-1} a_{\zeta}(k)
&=
\sum_{l=0}^{n/d-1} \sum_{k=0}^{d-1} a_{\zeta}(ld+k) \\
&=
\sum_{l=0}^{n/d-1} \sum_{k=0}^{d-1} \sum_{j=0}^{ld+k} c_q(j) c_q(ld+k-j) \\
&=
\sum_{l=0}^{n/d-1} \sum_{k=0}^{d-1} \sum_{i=0}^l \Big( c_{\zeta}(id) c_{\zeta}\big( (l-i)d \big) \Big) a_{\zeta}(k) \\
&=
\sum_{k=0}^{d-1} a_{\zeta}(k) \Big(\sum_{l=0}^{n/d-1}\sum_{i=0}^l c_{\zeta}(id) c_{\zeta}\big( (l-i)d \big) \Big) \\
&= 0.
\end{split}
\]
Thus $\sum_{k=0}^{n-1} a_q(k)$ is divisible by the cyclotomic polynomial $\Phi_d(q)$ for any divisor
$d>1$ of $n$. Hence $\sum_{k=0}^{n-1} a_q(k)$ is divisible by
\[
\prod_{1>d\mid n} \Phi_d(q) = [n],
\]
which completes the proof.
\section{Proof of Theorem~\ref{thm:main-2}}\label{sec:thm2-proof}
\noindent
Note first that using the formula~(\ref{Eq:GZ-1-1}) we have
\begin{equation}\label{Eq:main-2-1}
\sum_{k=0}^\infty a_q(k) = \fr{(-q^2;q^4)_{\infty}^2}{(1-q)^2(-q^4;q^4)_{\infty}^2}.
\end{equation}
The result is clear for $n=1$. Suppose that $n>1$ is odd. Let $d$ be an odd divisor of $n$ and let
$\zeta$ be a primitive root of unity of degree $d$.
In this case we have as was established in~\cite{Guo-Zudilin-2},
\begin{equation}\label{Eq:c-ld-2}
\lim_{q\to\zeta}c_q(ld) = c_{\zeta}(ld) = \fr{1}{4^l}{2l\choose l}.
\end{equation}
On the one hand, we have with the help of~(\ref{Eq:zeta}),
\[
\lim_{q\to\zeta} \fr{(-q^2;q^4)_{\infty}^2}{(1-q)(-q^4;q^4)_{\infty}^2}
=
\Big( \lim_{q\to\zeta} \fr{(-q^2;q^4)_{ld+k}}{(1-q)(-q^4;q^4)_{ld+k}} \Big)^2 \\
=
\fr{(-\zeta^2;\zeta^4)_k ^4}{(1-\zeta)^2 (-\zeta^4;\zeta^4)_k ^4}
\]
for any nonnegative $l$ and $0\leq k < d-1$. This shows that the right-hand-side
of~(\ref{Eq:main-2-1}) is bounded above by
\[
1+ \fr{1}{(1-\zeta)^2} \max_{0\leq k <d} \fr{(-\zeta^2;\zeta^4)_k ^4}{(1-\zeta)^2 (-\zeta^4;\zeta^4)_k ^4}.
\]
On the other hand, letting $q\to \zeta$ the left-hand-side now becomes
\[
\sum_{l=0}^{\infty} \Big(\sum_{i=0}^l c_{\zeta}(id) c_{\zeta} \big( (l-i)d \big) \Big) \sum_{k=0}^{d-1} a_{\zeta}(k)
\]
\[
=
\sum_{l=0}^{\infty} \Big(\sum_{i=0}^l \fr{1}{4^i}{2i\choose i}\fr{1}{4^(l-i)}{2(l-i)\choose l-i} \Big) \sum_{k=0}^{d-1} a_{\zeta}(k)
\]
\[
=
\left(\sum_{l=0}^{\infty} \fr{1}{4^l}{2l\choose l}+\sum_{l=0}^{\infty} \Big(\sum_{i=0}^{l-1} \fr{1}{4^i}{2i\choose i}\fr{1}{4^(l-i)}{2(l-i)\choose l-i} \Big)\right) \sum_{k=0}^{d-1} a_{\zeta}(k).
\]
As the series $\sum_{l=0}^{\infty} \fr{1}{4^l}{2l\choose l}$ diverges, we combine the foregoing identity with the fact that the left-hand-side is bounded
to deduce that
\[
\sum_{k=0}^{d-1} a_{\zeta}(k) = 0.
\]
Now we proceed as in the proof of Theorem~\ref{thm:main-1} to get the desired congruence.
\section{Proof of Theorem~\ref{thm:main-3}}\label{sec:thm3-proof}
By taking limits as $q\to \zeta$ on both sides of (\ref{Eq:GZ-3-1}) and following the same steps as in the proof of Theorem~\ref{thm:main-2}, we arrive at
\begin{equation*}
\sum_{k=0}^{n-1} a_{\zeta}(k) = 0
\end{equation*}
for any $d$-th primitive root of unity $\zeta$ of odd degree $1> d\mid n$. This yields
\begin{equation}\label{Eq:main-3-1}
\sum_{k=0}^{n-1} a_{q}(k) \equiv 0 \equiv q^{n-1} [n]^2 
\pmod{[n]}.
\end{equation}
Besides, by Guo~and~Zudilin~\cite[Lemma 3.1]{Guo-Zudilin-2}, the truncated sum
of~(\ref{Eq:GZ-3-1}) becomes
\[
\sum_{k=0}^{n-1} c_q(k)
=
\sum_{k=0}^{n-1} \fr{(aq;q^2)_k (q/a;q^2)_k (q;q^2)_{2k}}{(aq^6;q^6)_k (q^6/a;q^6)_k (q^2;q^2)_{2k}}
[8k+1]q^{2k^2}
= q^{\fr{n-1}{2}} [n] \Big(\fr{-3}{n}\Big)
\]
for $a=q^n$ or $a=q^{-n}$. Moreover, with these values of $a$ the coefficients satisfy
$c_q(k)=0$ whenever $\fr{n-1}{2}<k \leq n-1$, so by virtue of Lemma~\ref{lem:elementary}(a) we have
with these choices of $a$
\begin{equation}\label{Eq:main-3-2}
\sum_{k=0}^{n-1} a_{q}(k) = \Big(\sum_{k=0}^{n-1} c_q(k) \Big)^2.
\end{equation}
Now the desired supercongruence follows by~(\ref{Eq:main-3-1}) and~(\ref{Eq:main-3-2}) since
the polynomials $[n], 1-aq^n$, and $a-q^n$ are relatively prime.
\section{Proof of Lemma~\ref{lem:elementary}}\label{sec:lem-proof}
Part (a) is immediate from the assumption. As to part (b), we have
\begin{equation}\label{Eq:lem-2}
%\begin{split}
\sum_{j=0}^{ld+k}c(j)c(ld+k-j)
= \sum_{j=0}^{d-1} c(j)c(ld+k-j) + \sum_{j=d}^{2d-1} c(j)c(ld+k-j)
\end{equation}
\[
+\ldots +
\sum_{j=(l-1)d}^{ld-1} c(j)c(ld+k-j) + \sum_{j=ld}^{ld+k} c(j)c(ld+k-j) .
\]
We now handle the individual terms in the foregoing identity. Let $r_j$ be the remainder of the
division of $j$ by $d$ for $j=0,1,\ldots, ld+k$. Then for $a=0,1,\ldots, l-1$ and $ad\leq j \leq (a+1)d-1$,
we have $j = ad +r_j$ with $0\leq r_j<d$. Then
\begin{equation}\label{Eq:lem-3}
\begin{split}
\sum_{j=ad}^{(a+1)d-1} c(j)c(ld+k-j)
&=
\sum_{j: r_j\leq k} c(ad+r_j) c\big( (l-a)d + k-r_j \big) \\
& \qquad + \sum_{j: r_j > k} c(r_j) c \big( (l-a)d + k-r_j \big) \\
&=
\sum_{j= 0}^{k} c(ad+j) c\big( (l-a)d + k-j \big) \\
& \qquad + \sum_{j: r_j >k} c(r_j) c\big( (l-a)d + k- r_j \big) \\
&=
c(ad)c\big( (l-a)d \big) \Big(\sum_{j=0}^k \fr{c(ad+j)}{c(ad)} \fr{c\big( (l-a)d + k-j \big)}{c\big( (l-a)d \big)} \\
& \qquad + \sum_{j: r_j >k} \fr{c(r_j)}{c(ad)} \fr{c\big( (l-a)d + k-r_j \big)}{c\big( (l-a)d \big)}  \Big) \\
&=
c(ad)c\big( (l-a)d \big) \sum_{j=0}^k c(j) c(k-j)  \\
&\qquad + c(ad)c\big( (l-a)d \big) \sum_{j: r_j >k}
\fr{c(r_j)}{c(ad)} \fr{c\big( (l-a)d + k-r_j \big)}{c\big( (l-a)d \big)}
\end{split}
\end{equation}
We now claim that
\begin{equation}\label{Eq:lem-5}
\sum_{r_j =k+1}^{d-1} \fr{c(r_j)}{c(ad)} \fr{c\big( (l-a)d + k-r_j \big)}{c\big( (l-a)d \big)} = 0.
\end{equation}
Note first that
the claim is clear if $k>\fr{d-1}{2}$ since $c(r_j)=0$ for $\fr{d-1}{2}<r_j\leq d-1$ by assumption.
By the same assumption, it easy to see that the terms in the foregoing sum for which
 $r_j> \fr{d-1}{2} \geq k$ vanish. Now suppose that $k<r_j<\fr{d-1}{2}$. It follows that
 $r_j-k < \fr{d-1}{2}$ and so, $d+k-r_j > \fr{d-1}{2}$. Thus we get
\begin{equation*}\label{Eq:lem-5}
 \begin{split}
c(r_j) \fr{c\big( (l-a)d + k-r_j \big)}{c\big( (l-a)d \big)}
&=
\fr{c\big((l-a-1)d \big)}{c\big( (l-a)d \big)}
\fr{ c\big((l-a-1)d +d+k-j \big)}{c\big((l-a-1)d \big)} \\
&= \fr{c\big((l-a-1)d \big)}{c\big( (l-a)d \big)} c\big(d+k-r_j \big) \\
&= 0.
\end{split}
 \end{equation*}
 This proves the claim.
 Similarly, we have
\begin{equation}\label{Eq:lem-6}
\sum_{j=ld}^{ld+k} c(j) c(ld+k-j) = c(0)c(ld)\sum_{j=0}^k c(j)c(k-j).
\end{equation}
Now combine the relations (\ref{Eq:lem-2})--(\ref{Eq:lem-6})
to complete the proof.
\section{Conjectures}\label{sec:conjecture} 
We close this work by the following two conjectures which are related to Corollary~\ref{cor:main-1} and Corollary~\ref{cor:main-1} respectively. They both 
were suggested to us by an anonymous referee and they seem to be true by computational evidence.
\begin{conjecture}
For any odd prime number $p$, we have
\[
\sum_{k=0}^{p-1}\fr{(-1)^k}{8^{k}} \sum_{j=0}^k
{2j\choose j}{2k-2j\choose k-j}(6j+1)(6k-6j+1)
\equiv -\fr{p}{2} \pmod{p^2}.
\]
\end{conjecture}
\begin{conjecture}
For any odd prime number $p$, we have
\[
\sum_{k=0}^{p-1}\fr{1}{4^{k}} \sum_{j=0}^k
{2j\choose j}{2k-2j\choose k-j}(6j+1)(6k-6j+1)
\equiv p \pmod{p^2}.
\]
\end{conjecture}
%
%%%%%%%%%%%%%%%%%%%%%%%%%%%%%%%%%%%%%%%%%%%%%%%%%%%%%%%%%%%%%%%%%%%%%%%%%%%%%%
\noindent{\bf Acknowledgment.} The author is very grateful to the referees
for valuable comments and interesting suggestions which have improved the quality of the paper.
%%%%%%%%%%%%%%%%%%%%%%%%%%%%%%%%%%%%%%%%%%%%%%%%%%%%%%%%%%%%%%%%%%%%%%%%%%%%%

%
\end{document}